%% file: small.tex
\documentstyle[epic,eepic,epsf]{article}

\input mac.tex

\begin{document}

\title{RSOS models and 
       Jantzen-Seitz representations  
       of 
       Hecke algebras at roots of unity.
       }

\author{
Omar {\sc Foda}\thanks{
  Department of Mathematics, University  of Melbourne,
  Parkville, Victoria 3052, Australia.}, 
Bernard {\sc Leclerc}\thanks{
  D\'epartement de Math\'ematiques, 
  Universit\'e  de Caen, 
  BP 5186, 14032 Caen Cedex, France.}, 
Masato {\sc Okado}\thanks{
  Department of Mathematical Sciences, 
  Faculty of Engineering Science, 
  Osaka University, 
  Osaka 560, Japan.}, \\
Jean-Yves {\sc Thibon}\thanks{
  Institut Gaspard Monge, 
  Universit\'e de Marne-la-Vall\'ee,
  93166 Noisy-le-Grand Cedex, France.}  \ and \ 
Trevor A. {\sc Welsh}$^{\ast}$
  }

\date{ }

\maketitle

\begin{abstract}
A special family of partitions occurs in two apparently unrelated
contexts: the evaluation of 1-dimensional configuration sums of 
certain RSOS models, and the modular representation theory
of symmetric groups or their Hecke algebras $H_m$.
We provide an explanation of this coincidence by showing how
the irreducible $H_m$-modules which remain irreducible
under restriction to $H_{m-1}$ (Jantzen-Seitz modules)
can be determined from the decomposition  of a tensor product
of  representations of $\slchap_n$.
\end{abstract}

\setcounter{secnumdepth}{10}

\section{Introduction}

The solution of a class of ``restricted-solid-on-solid" (RSOS) models
by the corner transfer matrix method leads to the evaluation of
weighted sums of combinatorial objects called 
paths \cite{ABF}. 
The Kyoto group realized that these combinatorial sums are 
branching functions of the affine Lie algebra $\slchap_2$ \cite{DJMO}, and
was able to define similar models associated with other affine 
Lie algebras, in particular to $\slchap_n$ \cite{JMO}.

For the models associated with the cosets 
$(\slchap_n)_1 \times (\slchap_n)_1 / (\slchap_n)_2$,
a different description of the branching functions as generating
series of certain sets of partitions has been obtained 
in \cite{FOW}, and was used to derive fermionic expressions
for the configuration sums.

It turns out that exactly the same partitions arise in the 
modular representation theory of the symmetric groups: as
conjectured by Jantzen and Seitz \cite{JS} and established
recently by Kleshchev \cite{Kl1}, such partitions label the irreducible
representations  of a symmetric group $\SG_m$ over a field 
of characteristic $n$ which remain irreducible under restriction 
to $\SG_{m-1}$.

The aim of this Letter is to provide an explanation of this
seemingly mysterious coincidence.

The first point is to replace symmetric groups in characteristic $n$
by Hecke algebras of type $A$ over $\C$ at an $n$th root of unity.
Indeed, the representation theories of $\FF_n[\SG_m]$ and 
$H_m(\sqrt[n]{1})$ have many formal similarities, but the 
consideration of Hecke algebras removes the restriction of 
$n$ being a prime, which does not appear on the statistical mechanics side.

Moreover, a connection between the representation theory
of $H_m(\sqrt[n]{1})$ and the level 1 representations of
the quantum affine algebra $U_q(\slchap_n)$ has been pointed
out in \cite{LLT}.
Building on a conjecture of \cite{LLT}, recently proved by
Ariki and Grojnowski, we show that the Jantzen-Seitz type
problem for $H_m(\sqrt[n]{1})$ is equivalent to the decomposition
via crystal bases of tensor products of level 1 $\slchap_n$-modules.
Using the results of \cite{FOW}, we can then characterize the
Hecke algebra modules of Jantzen-Seitz type and explain 
the occurence of their partition labels in the configuration sums
of RSOS-models.  

As an application, we express the generating function 
of the Jantzen-Seitz partitions having a given
$n$-core in terms of branching functions of $\slchap_n$.
This is to be compared with a well-known result on
blocks of Hecke algebras.
Indeed, the blocks of $H_m(\sqrt[n]{1})$ are
labelled by $n$-cores, and the dimension
of a block is the number of $n$-regular partitions of $m$
with the corresponding $n$-core.
Using a formula first proved in the fifties by Robinson 
(in the symmetric group case) one can compute the generating series
of the dimensions of all blocks labelled by a given
$n$-core, and recognize the string function of the level 1
$\slchap_n$-modules.
Our result shows that some branching functions of $\slchap_n$ other
than the level 1 string function arise in a natural way in the 
representation theory of $H_m(\sqrt[n]{1})$.

\section{Characters and branching functions of $\slchap_n$}

Using the notion of paths, it was shown in \cite{DJKMO} that
the characters of the integrable highest weight modules
of $\slchap_n$ may be obtained by enumerating certain
{\it coloured multipartitions}.
In this Letter, we are interested only in the case when the
highest weight is of level one, and therefore the multipartitions
are simply partitions.

As is usual, we define a partition $\lambda$ of $m$ to be a
sequence 
$$\lambda=(\lambda_1,\lambda_2,\ldots,\lambda_k)$$
such that $\lambda_1\ge\lambda_2\ge\cdots\ge\lambda_k$
and $\lambda_1+\lambda_2+\cdots+\lambda_k=m$.
If $i>k$, we understand that $\lambda_i=0$.
The set of all partitions of $m$ is denoted $\Pi(m)$ and we write
$$\Pi=\bigcup_{m\ge0}\Pi(m).$$
Occasionally, it will be convenient to use the exponent notation
$$\lambda=(\lambda_1^{a_1},\lambda_2^{a_2},\ldots,\lambda_r^{a_r})$$
where here $\lambda_1>\lambda_2>\cdots>\lambda_r>0$, and $a_i>0$
specifies the multiplicity of $\lambda_i$ in $\lambda$.
If $a_i<n$ for $1\le i\le r$ then we say that the partition $\lambda$
is $n$-regular.
The set of all $n$-regular partitions of $m$ is denoted $\Pi_n(m)$,
and we define
$$\Pi_n=\bigcup_{m\ge0}\Pi_n(m).$$
The Young diagram associated with the partition $\lambda$ is
an array of $k$ left-adjusted rows of nodes (or boxes) in which
the $i$th row contains $\lambda_i$ nodes.
For example if $\lambda=(4,3,1)$, the corresponding Young diagram is:
$$
F^{(4,3,1)}=\smyoungd{
\multispan9\hrulefill\cr &&&&&&&&\cr
\multispan9\hrulefill\cr &&&&&&\cr
\multispan7\hrulefill\cr &&\cr
\multispan3\hrulefill\cr}.
$$
We will not distinguish between a partition and its Young diagram.
The partition $\lambda^\prime=(\lambda^\prime_1,\lambda^\prime_2,\ldots)$
conjugate to $\lambda$ is defined such that $\lambda^\prime_j$ is
the length of the $j$th column of $\lambda$, reading from left to right.

A coloured partition is a Young diagram in which each node
$\gamma$ is filled by its {\it colour charge} $c(\gamma)$
given by $c(\gamma)=(j-i)\mod n$, when $\gamma$ is the node
at the intersection of the $i$th row and the $j$th column.
For example, in the case $n=3$, the coloured partition
$\lambda=(5,5,4,1,1)$ appears as follows:
$$
\smyoungd{
\multispan{11}\hrulefill\cr &0&&1&&2&&0&&1&\cr
\multispan{11}\hrulefill\cr &2&&0&&1&&2&&0&\cr
\multispan{11}\hrulefill\cr &1&&2&&0&&1&\cr
\multispan{9}\hrulefill\cr &0&\cr
\multispan{3}\hrulefill\cr &2&\cr
\multispan{3}\hrulefill\cr }.
$$
Let $m_i$ be the number of nodes of $\lambda$ with colour charge $i$. 
The energy of $\lambda$ is defined by $E(\lambda)=m_0$,
and its weight by
$\wt(\lambda)=\La_0-\sum_{i=0}^{n-1}m_i\alpha_i$.
(We use freely the standard notations for roots, weights, etc. of 
$\slchap_n$, see {\it e.g.} \cite{DJKMO}.)
Then, 
the formal character of the basic representation $V(\Lambda_0)$ of
$\slchap_n$ is
$$
\ch V(\Lambda_0)=\sum_{\lambda\in\Pi_n}
e^{\wt(\lambda)}.
$$
A realisation of $V(\Lambda_0)$ that has a basis naturally indexed
by the set of $n$-regular partitions $\Pi_n$ was given in
\cite{DJKMO}.

In \cite{JMMO}, it was also shown how to describe combinatorially
the branching
functions of tensor products of irreducible $\slchap_n$-modules.
The branching function
$b_{\Lambda,\Lambda^\prime}^{\Lambda^{\prime\prime}}(q)$
is defined such that, for each $i\ge0$, the multiplicity of the
module $V(\Lambda^{\prime\prime}-i\delta)$
in the tensor product $V(\Lambda)\otimes V(\Lambda^\prime)$
is given by the coefficient of $q^i$ in
$b_{\Lambda,\Lambda^\prime}^{\Lambda^{\prime\prime}}(q)$.
Therefore, in terms of characters, we have
$$
\ch(V(\Lambda))\ch(V(\Lambda^\prime))=
\sum_{\Lambda^{\prime\prime}\in P^+}
b^{\Lambda^{\prime\prime}}_{\Lambda,\Lambda^\prime}(e^{-\delta})
\ch(V(\Lambda^{\prime\prime})).
$$
The result of \cite{JMMO} relies on the notion of a path
which we describe in the case where $\Lambda^\prime=\Lambda_0$.
A path $p$ is a sequence $p=(p_0,p_1,\ldots)$,
where $p_k\in P$, the weight lattice of $\slchap_n'$.
Given a dominant integral weight $\Lambda$ of level $l$ for some $l\ge0$,
and $\lambda\in\Pi_n$, a path $p=p(\lambda)$ is determined as follows.
For $k\ge\lambda_1$, define $p_k=\Lambda+\Lambda_k$.
Then for $0<k\le\lambda_1$, recursively obtain $p_{k-1}$ from $p_k$
by
$$
p_{k-1}=p_k-\epsilon_{k-1-\lambda^\prime_k},
$$
where $\epsilon_i = \Lambda_{i+1}-\Lambda_i$ and as usual the indices
are understood modulo $n$.
We write $\lambda \in {\cal Y}(\Lambda,\Lambda_0)$ 
if and only if all the coordinates $p_k$ of the path $p(\lambda)$
are dominant weights.
\begin{theorem}\label{BRCHtheorem}{\rm\cite{JMMO}}
$$
b^{\Lambda^\prime}_{\Lambda,\Lambda_0}(q)
=\sum_{
\scriptstyle\lambda \in{\cal Y}(\Lambda,\Lambda_0)\atop
\scriptstyle\wt(\lambda)=\Lambda^\prime-\Lambda\,(\mod\delta)}
q^{E(\lambda)}.
$$
\end{theorem}

\noindent
In the case when $\Lambda$ is of level one, the partitions
appearing in this sum were  characterised in \cite{FOW}.
We define the set $FOW(n,j)\subset\Pi_n$ as follows. Let
$\lambda=(\lambda_1^{a_1},\lambda_2^{a_2},\ldots,\lambda_r^{a_r})\in\Pi_n$.
Then $\lambda\in FOW(n,j)$ if and only if either $r=1$ or
$$
a_i+\lambda_i-\lambda_{i+1}+a_{i+1}=0\pmod n,
$$
for $i=1,2,\ldots,r-1$ and $j=(\lambda_1-a_1)\,\mod n$.

\begin{theorem}\label{FOWtheorem}{\rm\cite{FOW}}
$$
{\cal Y}(\Lambda_j,\Lambda_0)=FOW(n,j).
$$
\end{theorem}

\noindent If we now define $FOW(n,j,k)$ to be the subset of $FOW(n,j)$
comprising those partitions $\lambda$ for which
$\wt(\lambda)=\Lambda_k+\Lambda_{j-k}-\Lambda_j\,(\mod\delta)$,
then Theorem \ref{BRCHtheorem} immediately yields the following:

\begin{corollary}\label{FOWcorollary}{\rm\cite{FOW}}
Let $0\le j,k<n$. Then the branching function
$b^{\Lambda_k+\Lambda_{j-k}}_{\Lambda_j,\Lambda_0}(q)$
is given by:
$$
b^{\Lambda_k+\Lambda_{j-k}}_{\Lambda_j,\Lambda_0}(q)
=\sum_{\scriptstyle\lambda \in FOW(n,j,k)}
q^{E(\lambda)}.
$$
\end{corollary}

\noindent By means of this expression, the following was obtained:

\begin{theorem}{\rm\cite{FOW}}
Let $0\le j<n$ and $\La=\La_s+\La_t$ with $0\le s\le t<n$.
In addition, let $C$ be the Cartan matrix of $\goth{sl}_n$,
$e_i$ be the $(n-1)$-dimensional unit vector
$(0,\cd,0,\stackrel{i}{1},0\cd,0)^t$ and we set $e_n=0$.
Then
\[
b_{\La_j\,\La_0}^{\La}(q)=\sum_m\frac{q^{m^tC^{-1}m-m^tC^{-1}e_{s-t+n}+st/n}}
{\prod_{i=1}^{n-1}(q)_{m_i}},
\]
where the sum is taken over all $m\in(\Zn)^{\times(n-1)}$ satisfying
$t+\sum_{i=1}^{n-1}im_i=0\,(\mod\,n)$, and
$
(q)_k=(1-q)(1-q^2)\cd(1-q^k).
$
\end{theorem}

\section{Modular representations}

The set of partitions $FOW(n, j, k)$ has a definite meaning in
modular representation theory.
Indeed,  $\bigcup_{j, k} FOW(n,j, k)$ labels
the modular representations of $\SG_m$ in characteristic $n$
that remain irreducible under reduction to $\SG_{m-1}$ \cite{JS,Kl1}.

However, in the case of $FOW(n, j, k)$, $n$ can be any positive
integer, whereas in the case of $\SG_{m}$, it has to be a prime
number. This difference can be removed by working in the
context of Hecke algebras at an $n$th root of unity,
where the Jantzen-Seitz problem still makes sense.

The Hecke algebra $H_m(v)$ of type $A_{m-1}$, is the $\C(v)$-algebra
generated by the elements $T_1,\ldots,T_{m-1}$ subject to the relations
\begin{eqnarray*}
&&T_iT_{i+1}T_i=T_{i+1}T_iT_{i+1};\\
&&T_iT_j=T_jT_i \qquad \vert{i-j}\vert>1;\\
&&T_i^2=(v-1)T_i+v.
\end{eqnarray*}
In the case $v=1$, $H_m(v)$ may be identified with the group 
algebra of the symmetric group $\SG_m$, through identifying $T_i$ 
with the simple transposition $(i,i+1)\in\SG_m$. In fact, for 
generic values of $v$, $H_m(v)$ is isomorphic to the group 
algebra of the symmetric group $\SG_m$.
Hence,  it is semisimple, and its 
irreducible representations are parametrised by $\Pi(m)$.
A convenient realisation of the representation labelled by
$\lambda$ is the Specht module $S^\lambda$ in which the
entries of the representation matrices 
of the generators are elements of $\Z[v]$.

In the non-generic case when $v$ is a primitive $n$th root of unity,
$H_m(v)$ is no longer semisimple in general, and its
representations are not necessarily completely reducible.
The full set of irreducible representations
is indexed by $\Pi_n(m)$ (see \cite{DJ1}). The irreducible
module labelled by $\mu\in\Pi_n(m)$ is denoted by $D^\mu$.

Define $JS(n)$ to be the set of all partitions
that label the irreducible representations of $H_m(\sqrt[n]{1})$
which remain irreducible under reduction to $H_{m-1}(\sqrt[n]{1})$.
One of the aims of this work is to show that the partitions
in $JS(n)$ are defined by the same conditions as in the case of
symmetric groups, and to explain why
$$
JS(n) = \bigcup_{j, k} FOW(n,j, k).
$$

To tackle this problem, it is convenient to consider, as was done
in \cite{LLT,Ar2}, the Grothen\-dieck group $G_0(H_m(\sqrt[n]{1}))$
of the category of finitely generated $H_m(\sqrt[n]{1})$-modules.
The elements of $G_0(H_m(\sqrt[n]{1}))$ are classes $[M]$ of modules,
where $[M_1]=[M_2]$ if and only if the simple composition factors of $M_1$
and $M_2$ occur with identical multiplicities
(the order of the composition factors in the series is disregarded).
The sum is defined by $[M]+[N]=[M\oplus N]$.
It is known that this is a free abelian group
with basis the set $\{[D^\mu]\}$ of classes of irreducible
$H_m(\sqrt[n]{1})$-modules.
Then define ${\cal G}(n)$ as the direct sum
${\cal G}(n)=\bigoplus_m G_0(H_m(\sqrt[n]{1}))$.

As observed in \cite{LLT}, ${\cal G}(n)$ can be identified with the basic
representation $V(\Lambda_0)$ of $\slchap_n$, the action of the
Chevalley generators $e_i$, $f_i$ being given by the $i$-restriction
and $i$-induction operators, as defined in the fifties by G. de B.
Robinson in the case of symmetric groups. As an integrable highest
weight module of an affine algebra, ${\cal G}(n)$ has a canonical basis
in the sense of Lusztig and Kashiwara. It turns out that
this canonical basis coincides with the natural basis given by
the classes $[D(\mu)]$ of the irreducible $H_m(\sqrt[n]{1})$-modules.

To prove this, and to compute explicitly the
canonical basis, one considers the Fock space representation
${\cal F}$ of $\slchap_n$, which contains $V(\Lambda_0)$ as
its highest irreducible component:
$$
{\cal F} = V(\Lambda_0)\oplus V_{\rm low} \ .
$$
The standard basis $(v_\lambda)$ of ${\cal F}$ is labelled
by all $\lambda\in\Pi$ and one has an $\slchap_n$-homomorphism
$$
\matrix{ 
d: & {\cal F} & \longrightarrow & V(\Lambda_0)\simeq{\cal F}/V_{\rm low} \cr
   &v_\lambda & \longmapsto     & v_\lambda \ {\rm mod} \ V_{\rm low} 
}
$$
If ${\cal F}$ is identified with the Grothendieck ring of
all $H_m(v)$ for a generic $v$ by writing $[S^\lambda]=v_\lambda$,
then $d$ coincides with the decomposition map of modular
representation theory. In order to introduce the canonical
basis, one needs to $q$-deform the picture and to consider
the $q$-Fock space representation of $U_q(\slchap_n)$.

\section{The Fock space representation of $U_q(\slchap_n)$}

The $q$-Fock space $\cal F$ of $U_q(\slchap_n)$ 
has been described in \cite{Hay} using a $q$-analogue of the
Clifford algebra.
Another realization was given in
\cite{S,KMS} in terms of semi-infinite $q$-wedge products.
Let $(v_\lambda)$ denote the standard basis of weight
vectors of ${\cal F}$.
The Fock space ${\cal F}$ affords an integrable representation 
of $U_q(\slchap_n)$ whose decomposition into irreducible 
highest weight modules is given by
\begin{equation}\label{DECFOCK}
{\cal F} = \bigoplus_{k\ge 0} V(\Lambda_0-k\delta)^{\oplus p(k)} \,.
\end{equation}
In \cite{MM}, the lower crystal basis of ${\cal F}$ and its
crystal graph structure was described.
Let $A\subset \Q(q)$ denote the ring of rational functions without
a pole at $q=0$.
The lower crystal basis at $q=0$ of $\F$ is the pair
$(L,B)$ where $L$ is the lower crystal lattice given by
$$
L=\bigoplus_{\lambda\in\Pi} A\,v_\lambda\,,
$$
and $B$ is the basis of the $\Q$-vector space $L/qL$ given by
$$
B=\{v_\lambda\mod qL\mid\lambda\in\Pi\}\,.
$$
The action of Kashiwara's operators
$\tilde{f_i}$ on the element $v_\lambda\in B$
corresponds to adding to $\lambda$ a certain node
of colour $i$ which is called a good addable $i$-node.
Likewise the action of $\tilde{e_i}$ 
corresponds to the removal from $\lambda$ of a certain node
of colour $i$ which is called a good removable $i$-node.
As observed in \cite{LLT}, these nodes are precisely
those used by Kleshchev in his modular branching rule
for symmetric groups \cite{Kl3}, whence the terminology.

For each $\lambda\in\Pi$, the largest integer $k$ such that
$\tilde{e_i}^k(\lambda)\ne0$ ({\it resp} $\tilde{f_i}^k(\lambda)\ne0$)
is denoted by $\varepsilon_i(\lambda)$
({\it resp.} $\varphi_i(\lambda)$).
The crystal graph $\Gamma_n$ of $\cal F$ is
disconnected and reflects the decomposition (\ref{DECFOCK}). 
The connected component of the
empty partition is the crystal graph of $V(\La_0)$
and its vertices are labelled by $\Pi_n$.

We denote by $\{G(\mu),\mu \in \Pi_n\}$ the lower global basis
of $V(\Lambda_0)$. It is the $\Q[q,q^{-1}]$-basis of the
integral form $V_\Q(\Lambda_0)$ of $V(\Lambda_0)$ characterized by
$$
G(\mu) \equiv v_\mu(\mod qL),\qquad \overline{G(\mu)} = G(\mu).
$$
Here, for $v=xv_\emptyset\in V(\Lambda_0)$,
$\overline{v}$ is defined by
$\overline{v}=\overline{x}v_\emptyset$,
where $x\mapsto\overline{x}$ is
the $\Q(q)$-ring automorphism  of
$U_q(\slchap_n)$ given by
$$
\overline{q}=q^{-1},\qquad\overline{q^h} = q^{-h}\,,\qquad
\overline{e_i} = e_i,\qquad\overline{f_i} = f_i,
$$
for all $h\in P^\vee$ and $0\le i<n$.

The upper global crystal basis $\{G^{up}(\mu),\mu \in \Pi_n\}$
of $V(\Lambda_0)$ \cite{Ka2}
is the basis adjoint to
$\{G(\mu)\}$ with respect to the inner product $\<\cdot,\cdot\>$ defined
by
$$
\<v_\emptyset ,v_\emptyset\>=1\,\qquad
\<q^hv,v^\prime\>=\<v,q^hv^\prime\>\ ,\qquad
\<e_iv,v^\prime\>=\<v,f_iv^\prime\>\,
$$
for all $v,v^\prime\in V(\Lambda_0)$, $h\in P^\vee$ and $0\le i<n$.

From now on we fix $n\ge 2$ and write $H_m$ instead of
$H_m(\sqrt[n]{1})$.
It has been conjectured in \cite{LLT} and
proved in \cite{Ar2} that at $q=1$,
the upper global crystal basis of the $U_q(\slchap_n)$-module $V(\La_0)$
coincides in the identification $V(\Lambda_0) \simeq {\cal G}(n)$
with the basis $([D^\lambda])$ of ${\cal G}(n)$.
This was also proved by Grojnowski, using the results of \cite{G}.
This implies that

\begin{theorem}\label{GROTHtheorem}
If $\lambda\in\Pi_n(m)$ and
the coefficients $c_{\lambda\mu}(q)$ are defined by:
$$
\left( \sum_{i=0}^{n-1} e_i \right) G^{up}(\lambda)
= \sum_{\mu\in\Pi_n(m-1)} c_{\lambda\mu}(q) G^{up}(\mu),
$$
then,
$$
[D^\lambda\!\downarrow^{H_m}_{H_{m-1}}]
= \bigoplus_{\mu\in\Pi_n(m-1)} c_{\lambda\mu}(1) [D^\mu].
$$
\end{theorem}

\section{Jantzen-Seitz modules of Hecke algebras}

We are now in position to prove the  following result.

\begin{theorem}\label{JStheorem}
Let $\lambda=(\lambda_1^{a_1},\ldots,\lambda_r^{a_r})\in\Pi_n(m)$.
Then $D^\lambda\!\downarrow^{H_m}_{H_{m-1}}$ is irreducible,
i.e. $\lambda\in JS(n)$, 
if and only  if either
$r=1$ or $a_i+\lambda_i-\lambda_{i+1}+a_{i+1}\equiv 0 (\mod n)$
for $i=1,2,\ldots,r-1$.
\end{theorem}

\Proof
By Theorem~\ref{GROTHtheorem}, $D^\lambda\!\downarrow$ is 
irreducible if and only if 
$e_i\,G^{up}(\lambda) = G^{up}(\nu)$
for some $\nu$.
It is known
(see \cite{Ka3}, eqs. 5.3.8. -- 5.3.10.) that
\begin{equation}\label{magic}
e_i\,G^{up}(\lambda) = [\varepsilon_i(\lambda)]_q
G^{up}(\tilde e_i \lambda) 
+ \sum_{\mu\ne\lambda} E_{\lambda\mu}^i (q) \,G^{up}(\mu)
\end{equation}
where $E_{\lambda\mu}^i(q)$ is a Laurent polynomial invariant under
$q\mapsto q^{-1}$ such that 
$$
{E_{\lambda\mu}^i(q) } \in q^{2-\varepsilon_i(\lambda)}\Z[q]\,.
$$
If $\varepsilon_i(\lambda)=1$, then 
$[\varepsilon_i(\lambda)]_q=1$, and 
$E^i_{\lambda\nu}(q)=0$ since for all $k$ the coefficients of $q^k$
and $q^{-k}$ in $E^i_{\lambda\mu}$ are equal.
Therefore
$e_i\,G^{up}(\lambda) = G^{up}(\nu)$ where $\nu = \tilde e_i \lambda$.
On the other hand, if $e_i\,G^{up}(\lambda) = G^{up}(\nu)$
for some $\nu$,
then
(\ref{magic})  implies that
$\varepsilon_i(\lambda)=1$ so that $\nu = \tilde e_i \lambda$.

Hence we have obtained that $D^\lambda\!\!\downarrow$ is 
irreducible if and only if
$\varepsilon_i(\lambda)=1$ for some $i\in\{0,1,\ldots ,n-1\}$
and $\varepsilon_j(\lambda) = 0$ for $j\ne i$.

Using crystal graph theory (see {\it e.g.} Lemma 5.1 of \cite{JMMO}),
we see that the coefficient of $z^d$ in
the branching function $b_{\Lambda^\prime\,\Lambda_0}^\Lambda(z)$
is the number of vertices $b\in B$ for which
$\wt(b)=\Lambda-d\delta$ and $\varepsilon_j(b)\le\<\Lambda^\prime,h_j\>$
for $j=0,1,\ldots,n-1$.
Therefore, in the case $\Lambda^\prime=\Lambda_i$,
the coefficient of $z^d$ in $b_{\Lambda_i\,\Lambda_0}^\Lambda(z)$
is given by the number of vertices $b$ 
of the crystal graph of $V(\Lambda_0)$ such that
$$
{\rm wt}(b) = \Lambda - \Lambda_i - d\delta\ ,\qquad
\varepsilon_i(b)\le1\ ,\qquad
\varepsilon_j(b)=0 \quad (j\not = i).
$$


It follows from this discussion that the partitions $\lambda$
such that $D^\lambda\!\!\downarrow$ is irreducible can be identified
with the set of vertices $b$ of $B$
which contribute to the tensor product branching function
$b_{\Lambda_i\,\Lambda_0}^\Lambda(z)$ for some $i$ and some $\Lambda$.
Note that this branching function is non-zero only if
$\Lambda=\Lambda_k+\Lambda_{i-k}$ for some $k$.
Then, by \ref{FOWcorollary}, the Theorem is proved.

\begin{example}{\rm
The argument considered in the above proof may be illustrated by
reference to the $\emptyset$-connected component of the
crystal graph $\Gamma_3$ of $\slchap_3$ given
(up to weight 8) in Fig.~1.
Here, those partitions $\lambda$ for which $\lambda\in JS(n)$
have been highlighted with an asterisk. As in the above proof, 
these partitions correspond to the nodes $b$ in this crystal 
graph for which, for some $j$, $\varepsilon_j(b)\le1$ and 
$\varepsilon_i(b)=0$ for all $i\ne j$.}
\end{example}

\medskip
Now define the set $JS(n,\mu,d)$ to be the subset of $JS(n)$
comprising those partitions with $n$-core $\mu$ and $n$-weight $d$
(see \cite{JK,FLOTW} for definitions of $n$-core and $n$-weight).
Then define the generating series
$$
\chi_{n,\mu}(z)=\sum_{d\ge0} \#JS(n,\mu,d) z^d.
$$

\begin{theorem}\label{JScor}
\par
{\rm (i)} If $\lambda\in JS(n)$ then the $n$-core $\mu$ of $\lambda$
is a rectangular partition $\mu=(k^l)$ such that $k+l\le n$
(it is assumed here that if either $k=0$ or $l=0$ then
$(k^l)$ means the empty partition).

{\rm (ii)} If $k\ne0\ne l$ then
$$
\chi_{n,(k^l)}(z)= 
z^{-s}\,b_{\Lambda_{k-l}\,\Lambda_0}^{\Lambda_k + \Lambda_{-l}}(z) 
\,,
$$
where $s=\min(k,l)$.

{\rm (iii)} 
$$
\chi_{n,\emptyset}(z)=
\left(\,\sum_{k=0}^{n-1}
b_{\Lambda_k\,\Lambda_0}^{\Lambda_k + \Lambda_0}(z)\right)-(n-1)\,.
$$
\end{theorem}

\Proof In the tensor product $V(\Lambda_i)\otimes V(\Lambda_0)$,
all highest weights are of the form
$\Lambda_k+\Lambda_{-l}-e\delta$ with $k-l=i\mod n$
(for later convenience we use $-l$ and not $l$ here).
We take $0\le k,l<n$ here and can also assume that $k\le -l\mod n$
whereupon $k+l\le n$, and $l=0$ only if $k=0$.
By definition, the multiplicity of $V(\Lambda_k+\Lambda_{-l}-e\delta)$
is given by the coefficient of $z^e$ in
$b_{\Lambda_i,\Lambda_0}^{\Lambda_k+\Lambda_{-l}}(z)$ and hence,
by Corollary \ref{FOWcorollary}, by the number of partitions
$\lambda\in JS(n)$
for which $\wt(\lambda)=\Lambda_k+\Lambda_{-l}-\Lambda_i-e\delta$.
We claim that the $n$-core of such a $\lambda$ is the rectangular
partition $\mu=(k^l)$, for which we calculate
$\wt(\mu)=\Lambda_k+\Lambda_{-l}-\Lambda_{k-l}-s\delta=
\Lambda_k+\Lambda_{-l}-\Lambda_i-s\delta$,
where $s=\min(k,l)$ is the multiplicity of the colour charge 0 in $\mu$.
This follows from the fact that for every string of weights
$\Lambda,\Lambda-\delta,\Lambda-2\delta,\ldots,$ of the
$\slchap_n$-module $V(\Lambda_0)$ (where $\Lambda+\delta$ is not
a weight of $V(\Lambda_0)$), those partitions having these
weights have the same $n$-core, and this $n$-core has weight $\Lambda$.
Then since $\mu=(k^l)$ with $k+l\le n$ is manifestly an $n$-core,
it follows that it is the $n$-core of $\lambda$, thence proving
part (i).

Part (ii) follows immediately since in the case $k\ne0$ (so that
$l\ne0$), the partitions enumerated by the branching function
$b_{\Lambda_{k-l}\Lambda_0}^{\Lambda_k+\Lambda_{-l}}$
are precisely those elements $\lambda\in JS(n)$ having weight 
$\wt(\lambda)=\Lambda_k+\Lambda_{-l}-\Lambda_{k-l}-e\delta$,
for some $e$, and hence $n$-core $(k^l)$.

Finally, for $k=0$ and arbitrary $l$, each partition counted by
$b_{\Lambda_{-l}\Lambda_{0}}^{\Lambda_{0}+\Lambda_{-l}}$
has empty $n$-core, and hence contributes
to $\chi_{n,\emptyset}$.
However, the empty partition occurs for each $l$, hence an adjustment
of $n-1$ is needed after summing over all $l$.
No other partition $\lambda$ is repeated since, as indicated by
Theorem \ref{FOWtheorem},
the $b_{\Lambda_{-l}\Lambda_{0}}^{\Lambda_{0}+\Lambda_{-l}}$
to which it contributes is uniquely determined by
$-l\mod n=(\lambda_1-a_1)\mod n$.
(The summation over $-l$ is replaced by one over $k$ to give
the final result).

\newpage
\vfill
\centerline{\epsfbox{sl3_0.eps}}
\vfill

\newpage

\begin{example}{\rm
To illustrate this result, consider again the case $n=3$,
where we have the following branching functions (to three terms):
\begin{equation}
\begin{array}{lll}
b_{\Lambda_0,\Lambda_0}^{2\Lambda_0}&=&1+q^2+\cdots;\\[1mm]
b_{\Lambda_0,\Lambda_0}^{\Lambda_1+\Lambda_2}&=&q+2q^2+2q^3+\cdots;\\[1mm]
b_{\Lambda_1,\Lambda_0}^{\Lambda_1+\Lambda_0}&=&1+q+2q^2+\cdots;\\[1mm]
b_{\Lambda_1,\Lambda_0}^{2\Lambda_2}&=&q+q^2+2q^3+\cdots;\\[1mm]
b_{\Lambda_2,\Lambda_0}^{\Lambda_2+\Lambda_0}&=&1+q+2q^2+\cdots;\\[1mm]
b_{\Lambda_2,\Lambda_0}^{2\Lambda_1}&=&q+q^2+2q^3+\cdots.
\end{array}
\end{equation}
These are calculated using Corollary \ref{FOWcorollary} which leads to 
the enumeration of the nodes of Fig.~1. labelled by asterisks. The 
only rectangular $3$-cores are $\emptyset$, $(1)$, $(2)$ and $(1^2)$.
Using Theorem \ref{JScor}, we thus obtain:
\begin{equation}
\begin{array}{lll}
\chi_{n,\emptyset}&=&1+2q+5q^2+\cdots;\cr
\chi_{n,(1)}&=&1+2q+2q^2+\cdots;\cr
\chi_{n,(2)}&=&1+q+2q^2+\cdots;\cr
\chi_{n,(1^2)}&=&1+q+2q^2+\cdots.\cr
\end{array}
\end{equation}
\noindent These correspond to the following sets:
\begin{equation}
\begin{array}{lll}
JS(3,\emptyset,0)&=&\{\emptyset\};\cr
JS(3,\emptyset,1)&=&\{(3),(21)\};\cr
JS(3,\emptyset,2)&=&\{(6),(51),(3^2),(41^2),(321)\},\cr
JS(3,(1),0)&=&\{(1)\};\cr
JS(3,(1),1)&=&\{(4),(2^2)\};\cr
JS(3,(1),2)&=&\{(7),(43)\},\cr
JS(3,(2),0)&=&\{(2)\};\cr
JS(3,(2),1)&=&\{(5)\};\cr
JS(3,(2),2)&=&\{(8),(3^21^2)\},\cr
JS(3,(1^2),0)&=&\{(1^2)\};\cr
JS(3,(1^2),1)&=&\{(32)\};\cr
JS(3,(1^2),2)&=&\{(62),(44)\}.\cr
\end{array}
\end{equation}}
\end{example}

\section{Discussion}

We posed and solved the  Jantzen-Seitz problem for Hecke algebras of type 
$A$. The solution is obtained by mapping the problem to a problem in exactly 
solvable lattice models, namely that of characterising the space of states 
of a certain class of restricted solid-on-solid model in terms of the states 
of the corresponding (unrestricted) solid-on-solid models. The latter was 
solved in \cite{FOW}. The relationship between the two problems is based on 
the fact that both can be formulated in the same language of representation 
theory of $q$-affine algebra. The background to these two problems is 
discussed in greater detail in \cite{FLOTW}. 

In a forthcoming paper \cite{FLOTW-AK}, we plan to discuss the Jantzen-Seitz 
problem in the context of type-B Hecke algebras, and more generally, of  
Ariki-Koike (cyclotomic) Hecke algebras.

\bigskip

\footnotesize
{\noindent \it Aknowledgements --- }
We wish to thank Prof. Christine Bessenrodt for discussions that initiated 
this work, and Dr. Ole Warnaar for an earlier collaboration on which it is 
partly based. This work was done while B. Leclerc, M. Okado, and J.-Y. Thibon 
were visiting the Department of Mathematics, The University of Melbourne. 
These visits were made possible by financial support of the Australian 
Research Council (ARC). The work of O. Foda and T. A. Welsh is also supported 
by the ARC. 

{\noindent \it Note added --- } 
A forthcoming preprint, \cite{BO}, contains (among other things) an elementary 
purely combinatorial proof of part $(i)$ of \ref{JScor}.

\end{document}

%% file: mac.tex
\newtheorem{example}{Example}[section]

\newtheorem{theorem}[example]{Theorem}
\newtheorem{corollary}[example]{Corollary}

\def\cd{\cdots}

\def\La{\Lambda}

\def\wt{\mbox{\sl wt}\,}

\def\Z{{\bf Z}}
\def\Zn{\Z_{\ge0}}

\def\FF{{\ssym F}}

\def\Z{{\ssym Z}}
\def\Q{{\ssym Q}}
\def\C{{\ssym C}}
\def\Proof{\medskip\noindent {\it Proof: }}

\def\adots{\mathinner{\mkern2mu\raise1pt\hbox{.} 
\mkern3mu\raise4pt\hbox{.}\mkern1mu\raise7pt\hbox{.}}}
\def\<{\langle}
\def\>{\rangle}

\def\SG{{\goth S}}

\def\mod{{\rm\,mod\,}}

\def\ch{{\rm ch\,}}
\def\Sl{{\goth sl}}

\def\slchap{\widehat{\goth sl}}

\def\F{{\cal F}}

\def\F{{\cal F}}

\def\wt{{\rm wt\,}}

\def\F{{\cal F}}

\def\boxit#1#2{\setbox1=\hbox{\kern#1{#2}\kern#1}%

\dimen1=\ht1 \advance\dimen1 by #1 \dimen2=\dp1 \advance\dimen2 by #1
\setbox1=\hbox{\vrule height\dimen1 depth\dimen2\box1\vrule}%
\setbox1=\vbox{\hrule\box1\hrule}%
\advance\dimen1 by .4pt \ht1=\dimen1
\advance\dimen2 by .4pt \dp1=\dimen2 \box1\relax}

\font\tensym=msbm10
\font\sevensym=msbm7
\font\fivesym=msbm5
\newfam\ssymfam
\textfont\ssymfam=\tensym
\scriptfont\ssymfam=\sevensym
\scriptscriptfont\ssymfam=\fivesym
\def\ssym{\fam\ssymfam\tensym}

\font\tengoth=eufm10
\font\sevengoth=eufm7
\font\fivegoth=eufm5

\newfam\gothfam
\textfont\gothfam=\tengoth
\scriptfont\gothfam=\sevengoth
\scriptscriptfont\gothfam=\fivegoth
\def\goth{\fam\gothfam\tengoth}

\def\mathalign#1{\hbox to 0pt{\hss$\vcenter{\openup\jot
\tabskip=0pt plus1fil \halign to \hsize
{\hfil$\displaystyle ##$\tabskip=0pt&&$\displaystyle ##$\hfil
\tabskip=0pt plus1fil \cr #1}}$\hss}}

\def\makestrut#1#2{{\dimen12=#2
\divide\dimen12 by 4\dimen11=\dimen12\multiply\dimen11 by 3
\global\setbox#1=\hbox{\vrule height\dimen11 depth\dimen12 width0pt}}}

\newdimen\tadhdimen \newdimen\tabhdimen \newdimen\vdimen
\newdimen\smtadhdimen \newdimen\smtabhdimen
\newcount\splurt
\newbox\tadstrut \newbox\tabstrut
\newbox\smtadstrut \newbox\smtabstrut

\def\setyoungsize#1#2{          
  \tadhdimen=#1\tabhdimen=#1\advance\tabhdimen by -0.4truept%
  \vdimen=#2%
  \makestrut\tadstrut\vdimen
  \advance\vdimen by -0.4pt%
  \makestrut\tabstrut\vdimen}

\def\setsmyoungsize#1#2{        
  \smtadhdimen=#1\smtabhdimen=#1\advance\smtabhdimen by -0.4truept%
  \vdimen=#2%
  \makestrut\smtadstrut\vdimen
  \advance\vdimen by -0.4pt%
  \makestrut\smtabstrut\vdimen}

\setyoungsize{15pt}{15pt}
\setsmyoungsize{9pt}{9pt}

\def\youngt#1{%
  \vcenter{\offinterlineskip
  \halign{&\copy\tadstrut\hbox to \tadhdimen{\hss$##$\hss}\cr #1}}}
\def\youngd#1{%
  \vcenter{\offinterlineskip
  \halign{&\vrule##&\copy\tabstrut\hbox to \tabhdimen{\hss$##$\hss}\cr #1}}}

\def\smyoungt#1{{\vcenter{\offinterlineskip
  \halign{&\copy\smtadstrut
  \hbox to \smtadhdimen{\hss$\scriptstyle ##$\hss}\cr #1}}}}
\def\smyoungd#1{{\vcenter{\offinterlineskip
  \halign{&\vrule##&\copy\smtabstrut
  \hbox to \smtabhdimen{\hss$\scriptstyle ##$\hss}\cr #1}}}}

\def\hdashfill{\leaders\hbox to 3pt{%
\hfil\vrule width1.5pt height0.4pt depth0pt}\hfill}

\newcommand{\bi}{\begin{itemize}}
\newcommand{\ei}{  \end{itemize}}
\newcommand{\be}{\begin{equation}}
\newcommand{\ee}{  \end{equation}}
\newcommand{\bea}{\begin{eqnarray}}
\newcommand{\eea}{  \end{eqnarray}}


%% file: small.bbl
\newpage
\small
\begin{thebibliography}{99}

\bibitem{ABF}{\sc G.E.~Andrews, R.J.~Baxter} and {\sc P.J.~Forrester},
{\it Eight-vertex SOS model and generalized Rogers-Ramanujan-type 
     identities},
J. Stat. Phys. {\bf 35} (1984), 193--266.

\bibitem{Ar2}{\sc S. Ariki},
{\it On the decomposition numbers of the Hecke algebra of $G(m,1,n)$},
preprint, 1996.

\bibitem{BaxterBook}{\sc R.J. Baxter},
{\it Exactly solved models in statistical mechanics},
(1984) Academic Press. 

\bibitem{DJ1} {\sc R. Dipper} and {\sc G.D. James},
{\it Representations of Hecke algebras of general linear groups},
Proc. London Math. Soc. {\bf 52} (1986), 20--52.

\bibitem{DJKMO} {\sc E. Date, M. Jimbo, A. Kuniba, T. Miwa} and {\sc M. Okado},
{\it Paths, Maya diagrams, and representations of $\slchap(r,C)$},
Adv. Stud. Pure Math. {\bf 19} (1989), 149--191.

\bibitem{DJMO} {\sc E. Date, M. Jimbo, T. Miwa} and {\sc M. Okado},
{\it Automorphic properties of local height probabilities for
integrable solid-on-solid models},
Phys. Rev. {\bf B 35} (1987), 2105--2107.

\bibitem{FLOTW} {\sc O. Foda, B. Leclerc, M. Okado, J.-Y. Thibon}
                and {\sc T.A.~Welsh}, 
{\it  Modular representations of Hecke algebras and solvable lattice models},
submitted to
{\it Geometric Analysis and Lie Theory in Mathematics and Physics},
Lecture Notes Series of the Australian Mathematical Society.

\bibitem{FLOTW-AK} {\sc O. Foda, B. Leclerc, M. Okado, J.-Y.
Thibon} and {\sc T.A.~Welsh}, in preparation.

\bibitem{FOW}{\sc O. Foda, M. Okado} and {\sc S.O. Warnaar},
{\it A proof of polynomial identities of type 
$\widehat{sl}(n)_1\otimes\widehat{sl}(n)_1/\widehat{sl}(n)_2$},
J. Math. Phys. {\bf 37} (1996), 965--986.

\bibitem{G}{\sc I. Grojnowski}, {\it Affine Hecke algebras 
(and affine quantum $GL_n$) at roots of unity},
Internat. Math. Res. Notices (1994), 215-217.

\bibitem{Hay}{\sc T. Hayashi},
{\it $q$-analogues of Clifford and Weyl algebras - spinor and
oscillator representations of quantum enveloping algebras},
Commun. Math. Phys. {\bf 127} (1990), 129--144.

\bibitem{Ja}{\sc G.D.~James}, {\it The decomposition
matrices of $GL_n(q)$ for $n\le 10$}, Proc. London
Math. Soc. {\bf 60} (1990), 225--265.

\bibitem{JK}{\sc G.D.~James} and {\sc A.~Kerber},
{\it The representation theory of the symmetric group},
Addison-Wesley, 1981.

\bibitem{JS}{\sc J.C. Jantzen} and {\sc G.M. Seitz},
{\it On the representation theory of the symmetric groups},
Proc. London Math. Soc. {\bf 65} (1992), 475--504.

\bibitem{JMMO}{\sc M. Jimbo, K. Misra, T. Miwa} and {\sc M. Okado},
{\it Combinatorics of representations of $U_q(\slchap(n))$ at $q=0$},
Commun. Math. Phys. {\bf 136} (1991), 543--566.

\bibitem{JMO}{\sc M. Jimbo, T. Miwa} and {\sc M. Okado},
{\it Solvable lattice models whose states are dominant integral
weights of $A_{n-1}$}
Lett. Math. Phys. {\bf 14} (1987) 123--131.

\bibitem{Ka2}{\sc M. Kashiwara},
{\it On crystal bases of the $q$-analogue of universal enveloping algebras},
Duke Math. J. {\bf 63} (1991), 465--516.

\bibitem{Ka3}{\sc M. Kashiwara},
{\it Global crystal bases of quantum groups},
Duke Math. J. {\bf 69} (1993), 455--485.

\bibitem{KMS}{\sc M. Kashiwara, T. Miwa} and {\sc E. Stern},
{\it Decomposition of $q$-deformed Fock spaces}, Selecta Mathematica 1996.

\bibitem{Kl1}{\sc A.S. Kleshchev},
{\it On restrictions of irreducible modular representations
of semisimple algebraic groups and symmetric groups to natural
subgroups I},
Proc. London Math. Soc. {\bf 69} (1994), 515-540.

\bibitem{Kl3}{\sc A.S. Kleshchev},
{\it Branching rules for the modular representations of symmetric groups III; 
some corollaries and a problem of Mullineux},
J. London Math. Soc. {\bf 2} (1995).

\bibitem{LLT}{\sc A. Lascoux, B. Leclerc} and {\sc J.-Y. Thibon},
{\it Hecke algebras at roots of unity and crystal bases of
quantum affine algebras},
Commun. Math. Phys. {\bf 181} (1996), 205-263.

\bibitem{MM}{\sc K.C. Misra} and {\sc T. Miwa},
{\it Crystal base of the basic representation of $U_q(\widehat{\Sl}_n)$},
Commun. Math. Phys. {\bf 134} (1990), 79--88.

\bibitem{S}{\sc E. Stern}, {\it Semi-infinite wedges and vertex operators},
Internat. Math. Res. Notices 1995, 201-220.

\bibitem{BO}{\sc Ch. Bessenrodt} and {\sc J. Olsson}, 
{\it Residue symbols and Jantzen-Seitz partitions}, 
to appear.

\end{thebibliography}
